\newcommand{\N}{\mathbb{N}}
\newcommand{\C}{\mathbb{C}}
\newcommand{\R}{\mathbb{R}}

\newcommand{\Z}{\mathbb{Z}}

\newcommand\cA{\mathcal{A}}
\newcommand\cB{\mathcal{B}}
\newcommand\cC{\mathcal{C}}

\newcommand\cH{\mathcal{H}}

\newcommand\cL{\mathcal{L}}
\newcommand\cM{\mathcal{M}}
\newcommand\cN{\mathcal{N}}

\newcommand\cP{\mathcal{P}}
\newcommand\cQ{\mathcal{Q}}
\newcommand\cR{\mathcal{R}}
\newcommand\cS{\mathcal{S}}

\newcommand\gG{\Gamma}

\newcommand\gs{\sigma}

\newcommand\rk{\operatorname{rk}} 
\newcommand\K{\operatorname{K}}

\newcommand{\Id}{\mathrm{Id\,}}

\newcommand{\ga}{\gamma}
\newcommand{\si}{\sigma}
\newcommand{\la}{\lambda}
\newcommand{\Ker}{\mathrm{Ker\,}}

\newcommand{\Ind}{\mathrm{Ind\,}}
\newcommand{\J}{J\,}

\newcommand{\cf}{\mathcal{F}}

\newcommand{\cg}{\mathcal{G}}
\newcommand{\cLL}{\mathcal{L_\lambda}}
\newcommand{\cBL}{\mathcal{B_\lambda}}
\newcommand{\hkr}{H^{k+s}(\Omega;\R^m )}
\newcommand{\hk}{H^{k+s}(\Omega;\C^m )}
\newcommand{\hs}{H^{s}(\Omega;\C^m )}

\newcommand{\hkman}[1]{H^{k+s}(#1;\C^m )}
\newcommand{\hsman}[1]{H^{s}(#1;\C^m )}

\newcommand{\hsbr}{H^{s}(\Omega;\R^m)\times H^{s*}(\partial\Omega;\R^r )}
\newcommand{\hbb}{H^{s*}(\partial \Omega;\C^r )}
\newcommand{\hbbr}{H^{s*}(\partial \Omega;\R^r )}

\newcommand{\diag}{\operatorname{diag}}

\newcommand{\ind}{\operatorname{ind}}

\newcommand{\coker}{\operatorname{coker}}

\newcommand{\hfor}{\hbox{ \,for\, }}

\newcommand{\hif}{\hbox{ \,if\,} }
\newcommand{\hof}{\hbox{\,of\,} }

\newcommand{\hforall}{\hbox{ \,for\  all\, }}

\newcommand{\sk}{\vskip10pt}

\newcommand{\ch}{\operatorname{ch}}

\newcommand{\bt}{\begin{theorem}}
\newcommand{\et}{\end{theorem}}
\newcommand{\bc}{\begin{corollary}}
\newcommand{\ec}{\end{corollary}}
\newcommand{\bp}{\begin{proposition}}
\newcommand{\ep}{\end{proposition}}
\newcommand{\bl}{\begin{lemma}}
\newcommand{\el}{\end{lemma}}
\newcommand{\br}{\begin{remark}}
\newcommand{\er}{\end{remark}}
\newcommand{\bd}{\begin{definition}}
\newcommand{\ed}{\end{definition}}
\newcommand{\be}{\begin{equation}}
\newcommand{\ee}{\end{equation}}
\newcommand{\ba}{\begin{array}}
\newcommand{\ea}{\end{array}}
\newcommand{\ra}{\rightarrow}
\newcommand{\eqr}{\eqref}
 
\newcommand{\cipe}{\cite{[Pe]}\ }
\documentclass[a4paper]{amsart}
\usepackage[mathscr]{eucal}

\usepackage{mathtools}
\usepackage{amsmath}
\usepackage{amsfonts}
\usepackage{amssymb}
\usepackage[applemac]{inputenc}
\input{diagxy}

\begin{document}
%\nocite{*}
\newtheorem{theorem}{Theorem}[section]
\newtheorem{corollary}[theorem]{Corollary}
\newtheorem{lemma}[theorem]{Lemma}
\newtheorem{proposition}[theorem]{Proposition}
\newtheorem{remark}{Remark}[section]
\newtheorem{definition}{Definition}[section]
\numberwithin{equation}{section}
\numberwithin{theorem}{subsection}
\numberwithin{definition}{subsection}
\numberwithin{remark}{subsection}
%\shorttitle{Index and bifurcation}

\title[Family index theorem and bifurcation]{The family index theorem and bifurcation of solutions of nonlinear elliptic BVP %boundary value problems
}
\author{J. Pejsachowicz }

\address{Dipartimento di Matematica Politecnico di Torino\\
    Corso Duca degli Abruzzi  24, 
       10129 Torino, Italy.}
\email{jacobo.pejsachowicz@polito.it}



%\date{\today}

\thanks{This work was supported by MIUR-PRIN 2007-Metodi variazionali e topologici nei fenomeni nonlineari}

%\tableofcontents 
\date{today}
\begin{abstract} 
We  obtain some  new  bifurcation criteria for  solutions of general  boundary value problems  for nonlinear elliptic systems of  partial differential equations.  The results are of different nature  from the ones that can be obtained via  the traditional Lyapunov-Schmidt reduction.  Our sufficient conditions for bifurcation  are derived from the Atiyah-Singer family index theorem  and therefore they depend only on the coefficients of derivatives of leading order of  the linearized differential operators. They  are computed explicitly from the coefficients without the need of solving the linearized  equations. Moreover,  they are stable  under  lower order perturbations.
\end{abstract}
\maketitle
\section{Introduction}\label{sec:0}
 
Let us begin with a rough description of the general setting and the background of the problem. 

 We consider a system  of nonlinear partial differential equations of  the form
 \begin{equation}\label{bvp1}
\left\{\begin{array}{l} \cf^i\,(\la, x,u,\ldots ,D^{k}u)=0 \hfor  x \in  \Omega , \, 1 \le  i \le m \\ \cg^i(\la, x,u,\ldots, D^{k_i} u)=0 \hfor  x \in  \partial \Omega , \, 1 \le  i \le r.\end{array}\right. \end{equation} 
Here $\Omega$ is  an open bounded subset of $\R^n$ with smooth boundary $\partial \Omega,$ $u\colon\bar\Omega\rightarrow \R^m$ is a vector function and $\la\in \R^{q}$ is a parameter.   $\cf^i, \cg^i $ are smooth functions defined on jets of order $k$  and $k_i$ respectively such that $\cf^i (\la,x, 0) = 0,\  \cg^i( \la,x,0)=0.$  We will not assume  that $k_i\leq k-1$ since we won't need this here. 

Because of the last assumption,  the function $u\equiv 0$ is a solution of problem  \eqref{bvp1} for every  $\la \in \R^q.$  The set $\ \R^q \times \{0\} $ is called a trivial branch of solutions  of \eqref{bvp1}. 
Nontrivial  solutions  of \eqref{bvp1} are solutions $(\la,u)$ with  $u\neq 0.$ Roughly speaking, a  {\it  bifurcation point  from the trivial branch} for solutions  of \eqref{bvp1} is a point   $\la \in \R^q$ such that arbitrarily  close to $(\la,0)$ there are nontrivial solutions of the equation. 
 
  For each $\la\in \R^q,$ the functions $\cf^i (\la,x, 0),\  \cg^i( \la,x,0) $ define a nonlinear differential operator  \begin{equation}\label{flambda}  (\cf_\la,\  \cg_\la) \colon C^\infty(\bar \Omega;R^m) \ra C^\infty(\bar \Omega;R^m)\times C^\infty(\partial\Omega;R^r).\end{equation}
  
  Let $\left(\cLL, \cBL\right)$ be the linearization of $(\cf_\la,\cg_\la)$ at $u\equiv 0.$  It is easy to see that a necessary condition for a point $ \la$ to be a bifurcation point  is  the existence of a non vanishing solution $v$  of the  linearized  problem  
\begin{equation}
\label{linsol}\left\{
\begin{array}{l}
 {\mathcal L}_\la(x,D) v(x) = 0,  \hfor  x \in  \Omega\\
\mathcal {B}_\la(x,D) v(x) =0, \hfor  x \in\partial \Omega.\end{array}\right.\end{equation} 

However, the above condition is not sufficient  and  the  goal  of the  linearized bifurcation theory is to obtain sufficient conditions  for the appearance of nontrivial solutions from invariants  associated to  the  linearization  $(\cLL, \cB_\la).$ 
  
   If  the  operator  $ \cLL$  is elliptic and the boundary operator  $ \cB_\la  $ verifies the Shapiro-Lopatinskij  condition with respect to $\cLL ,$   then for all $\la \in \R^q,$   the differential operator $\left(\cLL , \cB_\la \right)$  induces a Fredholm operator between  the Hardy-Sobolev spaces naturally associated to the problem.  It follows from this  that  $\cf^i,\cg^i$   define on  a neighborhood of zero  a  family of nonlinear Fredholm maps. The Fredholm property is essential. It gives the possibility  to recast, at least locally, the study of a bifurcation problem  to an equivalent problem for a finite number of nonlinear equations in a finite number of indeterminates. This is  the essence of the celebrated Lyapunov-Schmidt method.  A  typical further assumption in this setting is that points where \eqref{linsol} holds  are isolated.  Assuming this, a number of sufficient conditions for bifurcation can  be obtained using  either  analytical or topological methods \cite{[Ki],[Ra-St],[Zv-Ra], [Iz], [Iz-1],[Rab]}. 
 
 In the past years  we  worked out  a different approach to bifurcation,  based  on various homotopy  invariants of families of Fredholm operators defined by the linearization along  the trivial branch.  The invariants under consideration, the index bundle, the spectral flow and others are borrowed from elliptic topology. They arise in bifurcation theory as a tool linking  the nontrivial topology of the parameter space with the appearance of nontrivial solutions of  the equation. 
 
In  \cite{[Pe]}  we introduced  an {\it index of bifurcation} which counts algebraically the bifurcation points  of a family of nonlinear Fredholm maps parametrized by  open subsets of a  compact manifold or polyhedron $\Lambda.$     The index takes  values in a finite abelian group $J(\Lambda)$ associated to the parameter space.   It has similar properties  to the well known fixed-point index. Namely, it possesses  the existence, additivity, excision and homotopy invariance property. At an isolated singular point of the linearization the index of bifurcation  can be computed using the  Lyapunov-Schmidt reduction. But, what is most important, the total index is derived from a well known elliptic invariant;  the index bundle  of the linearization along the trivial branch.   The Atiyah-Singer family index theorem allows us  to compute the index bundle directly from  the principal part of the linearization, i.e.,  the coefficients of  the leading  derivatives of   $\cLL(x,D) $ and $\cB_\la (x,D).$

 For families of elliptic boundary value problems parametrized by $\R^q,$  viewed as an open subspace of the sphere $S^q,$ the results are particularly simple. The groups $J(S^q)$ are finite  cyclic groups whose orders  have been computed by Adams and others.   In  Theorem 1.4.1 of \cite{[Pe]},  under the  assumption   that  the principal part of the boundary operator  $\cB_\la(x,D)$ is  independent of $\la$ and that the principal part of the  interior operator $\cLL(x,D)$ is independent of  $\la$ near the boundary, we computed the index of bifurcation obtaining in this way sufficient conditions for bifurcation in terms of a number defined  as  the integral  of a differential form constructed explicitly from the principal symbol of the linearization $\cLL(x,D).$

 The purpose of the present paper  is to extend  Theorem 1.4.1 of \cite{[Pe]}  to the case  in which also the principal symbol of the  boundary operator is  parameter dependent.  
 
In \cite {[Pe]} we proved  a parametrized version of  the Agranovich reduction  in order to recast the calculation of the index bundle of the  family of linearizations along the trivial branch   to that  of a family of pseudo-differential operators on $\R^n.$ Then we  applied  the  Atiyah-Singer family index theorem to the latter. In this article,  in addition to our previous results,  we will need yet another type of reduction. We will show that  the index bundle of a  family  of elliptic  boundary value problems  whose interior operator is independent from the parameters coincides with  that of  a family of pseudo-differential operators on the boundary.  In the comparison of  two boundary value problems with the same interior operator  the latter is known under the name of  Agranovich-Dynin reduction \cite{[Ag-Dy]},\cite{[R-S]}. 

Our result will be  proved by combining both reductions with the cohomological form of the family index theorem.   Consequently, our criteria for the existence of  bifurcation points will be formulated  in terms of  the Bott-Fedosov  degree of  two maps  $\gs$  and  $\tau,$ with values in $GL(m;\C)$ and $GL(r;\C)$ respectively,  which  are naturally  attached to the reductions discussed above. 

 The map $\gs$ is constructed out of the principal symbol of the family of  interior  operators $\cL$ while $\tau$ is constructed by restricting the principal symbol of the boundary operators to the  vector bundle $M_+$ of all  stable solutions of a a family of ordinary  differential equations  canonically associated to $\cL.$

Extending the operators $\cLL(x,D)$ to the double of $\bar \Omega$ would give us  sufficient conditions for bifurcation in a slightly different but equivalent form.  We consider the chosen approach more straightforward.   Notice that we still have to assume independence from the parameter  $\la$ of the coefficients of the  operators $\cLL(x,D)$  near to  the boundary of $\Omega.$  Taking into account the variation on the boundary would lead us to pseudo-differential operators with operator valued symbols.

The paper is organized as follows. In Section 2 we state our main result. Section 3 is a short review of  the background  material from \cite{[Pe]}. Section 4 contains  the  proofs and in Section 5 we construct some examples  of semi-linear elliptic boundary value problems illustrating  our bifurcation result. An appendix is devoted to the comparison of the Fedosov's approach to Bott-Fedosov degree in \cite{[Fe]} with the one used by Atiyah and Singer  in \cite{[At-Si-III]}.   

 \sk 
 \section{Statement of the main theorem}
Let   
 $(\cf,\cg) \colon \R^{q}\times C^{\infty}(\bar\Omega; \R^{m}) \to  C^\infty(\bar \Omega;R^m)\times C^\infty(\partial\Omega;R^r)$ be a parametrized family of nonlinear differential operators  defined by
  \begin{equation} \label {fg}\left\{ \displaystyle \begin{array}{lll} \cf(\la, u) &=&\left( \cf^1(\la,x,u,\dots, D^k u), \dots,\cf^m( \la,x,u,\dots, D^k u)\right)\\
 \cg(\la, u) &=&\left ( \cg^1(\la,x,u,\dots, D^{k_1} u), \dots,\cg^r( \la,x,u,\dots, D^{k_r} u)\right).\end{array}\right.
 \end{equation} 
 
 Keeping our notations from \cite{[Pe]} we  will denote with $(\cf_\la,\cg_\la)$ the operator corresponding to the parameter value $\la \in \R^q.$ In general families of  differential operators will be denoted  using  calligraphic letters while the families of induced operators on Hardy-Sobolev spaces will be denoted with the corresponding roman capitals.  For example, in our notations $\cLL=\cLL(x,D);$  the induced operator being $L_\la.$ 
 \sk
Denoting by  $v_{j\alpha}$ the variable corresponding to $D^\alpha u_j,$ for each  fixed $\la,$ {\it the  linearization} of  $(\cf_\la,\cg_\la)$  at  $u\equiv 0$  is the linear  differential operator $(\cLL,\cB_\la)$ defined by 
\begin{equation}
\label{lin}\left\{\displaystyle
\begin{array}{lll} {\cL}_\la (x,D) &=&\sum_{|\alpha| \leq k} a_{\alpha}(\la, x) D^{\alpha}, \\
{\cB}^i_{\la}(x,D)&=&\ga_0 \sum_{ |\alpha| \leq k_i} b_\alpha (\la,x)D^{\alpha},\,1\leq i \leq r, 
\end{array}\right.
\end{equation}
where  the  $ij$-entries of
 $a_{\alpha } \in C^\infty( \Lambda\times \bar\Omega; \R^{m\times m} )$ and  $b^i_{\alpha } \in C^\infty( \Lambda\times \bar\Omega; \R^{1\times m} )$ are
\begin{equation}
\label{lin1}
 a^{ij}_{\alpha }(\la,x) ={ \frac{\partial \cf^i}{\partial v_{j\alpha}}}(\la,x,0), \    b^{ij}_{\alpha }(\la,x) ={ \frac{\partial \cg^i}{\partial v_{j\alpha}}}(\la,x,0),
\end{equation}
and  $\ga_0$ is  the operator "restriction to the boundary".

\sk
We assume:

\begin{itemize}
 \item[$H_1)$]  For all $\la \in R^q,$  the linearization $(\mathcal L_\la, \mathcal B_\la)$  of $(\cf_\la, \cg_\la)$  at $u\equiv 0,$
defines  an elliptic boundary value problem. Namely, $\cLL$ is elliptic, properly elliptic at the boundary and $\mathcal B_\la$  verifies the Shapiro-Lopatinskij condition with respect to $\cLL$ (see \cite[Definition 5.2.1]{[Pe]}). 
\smallskip
 \item [$H_2)$] 
 The coefficients   $a^{ij}_\alpha, b^{ij}_\alpha $ of the linearized family  $(\mathcal L, \mathcal B)$   extend  to  smooth functions  defined on  $S^q \times \bar \Omega,$ where  $S^q=\R^q\cup\{\infty\} $ is the one point compactification of $\R^q.$ Moreover
the problem 
  \[ \left\{ \begin{array}{l} {\mathcal L}_\infty (x,D) u(x)=\displaystyle \sum_{|\alpha|\leq k} a_{\alpha }(\infty, x) D^{\alpha }u(x) =f(x), \  x \in \Omega\\
\mathcal {B}^i_\infty (x,D) u(x)=\displaystyle \ga_0\sum_{|\alpha|\leq k_i} b^i_\alpha (\infty,x)D^{\alpha }u(x)=g(x),\  x\in \partial\Omega, \, 1\leq i\leq r, \end{array} \right .\]  is elliptic and has a 
unique solution for every $f \in C^\infty(\bar \Omega; \R^m)$ and every  $g\in C^\infty (\partial \Omega; \R^r ).$

\smallskip 
   \item [$H_3)$]  
 The restrictions   of  the coefficients   of the  leading terms  of  ${\mathcal L}_\la (x,D)$ to a neighborhood  of $ \partial \Omega$  are independent of $\la.$ Moreover the principal symbol of the operator $\cL_\infty(x,D)$ commutes with the principal symbol of $\cLL(x,D)$  for all $\la\in S^q.$ 
  \end{itemize} 
  
  \begin{remark} {\rm The hypothesis  $H_2, H_3 $ are  restrictive.  We will discuss  elsewhere  bifurcation of  elliptic boundary value problems  under  different assumptions, which do not require the extension of the linearized family to $S^q.$ Obviously, the  principal  symbols of $\cL_\infty$  and $\cLL $   commute   if either the principal symbol of $\cL$ is constant or the principal symbol of  $\cL_\infty $ is diagonal.  This later condition is not needed if the principal symbol of  $\cB$ is independent from the parameter  \cite{[Pe]}.}\end{remark}
  
 \begin{definition} \label{def:1}
  A bifurcation point from the trivial branch  for solutions of \eqref{bvp1} is a point $ \la_*$ such that there exists    a sequence  $(\la_n, u_n)\in \R^q \times C^\infty (\bar\Omega;\R^m)$ of  nontrivial solutions  of  \eqref{bvp1}  with   $\la_n\rightarrow \la_*$ and $u_n\rightarrow 0 $ uniformly together with all of their derivatives. 
\end{definition}

\sk
 
 A sufficient condition  for the existence of bifurcation points of \eqref{bvp1}  is  that  index of bifurcation of the family of Fredholm maps between Hardy-Sobolev spaces induced by $(\cf,\cg)$ does not vanish (see \cite{[Pe]}). In Section $4$ we will compute the  index  from two integers  associated to the linearization $(\cL,\cB)$  at the trivial branch, called the interior and  the boundary multiplicity, together with some natural numbers $n(q)$ related to the order of the group $\J(S^q).$  As we said in the introduction,  the interior and boundary multiplicity  will be defined as  the Bott-Fedosov degree of  two maps  $\gs$  and  $\tau$ with values in $GL(m;\C)$ and $GL(r;\C)$ respectively. 
\sk
The construction of $\gs$ is as follows:

 Let $p(\la,x,\xi ) \equiv \sum^{}_{\vert \alpha \vert =k} a_{\alpha }(\la,x)\xi ^{\alpha }$ be the principal symbol of $\cLL.$  Since  the principal symbol  is obtained  substituting  the operator  $D_j= -i \frac{\partial }{\partial x_j}$ with
 the variable $\xi_j,$  $p(\la,x,\xi )$ is  a complex matrix verifying  the reality condition \begin{equation} \label{real} p\,(\la,x,-\xi)=\bar p\,(\la,x,\xi).\end{equation}  
 
  By ellipticity, $ p(\la,x,\xi )\in GL(m;\C)$ if $\xi\neq 0.$ On the other hand, by $H_3,$
   $ p (\la,x,\xi )= p (\infty,x,\xi )$   for all $x$ in a neighborhood of $\partial\Omega.$  Hence putting  
   \[\sigma(\la,x,\xi)= \begin{cases} p(\la,x,\xi) p(\infty,x,\xi)^{-1}& \hif\,   x\in \Omega, \,\xi \neq 0  \\Id  & \hif \, x \notin \Omega,\end{cases} \]
  we get a smooth map 
 
  \begin{equation} 
  \label{symbo}   
\sigma \colon S^q\times (\R^{2n} - \Omega \times  \{0\})  \to GL(m;\C).\end{equation} 

\sk

Now let us define  the matrix function  $\tau.$

 We take a neighborhood $\cN$ of $\partial\Omega$  of the form $\cN \simeq \partial\Omega\times (-1,1).$   At a  point $x$ belonging to $\Omega\cap\cN$ we  will use a coordinate system of the  form  $(x',t),$ where  $x' =(x'_1,\dots,x'_{n-1})$ is a coordinate  system on  the manifold $\gG =\partial \Omega$  and $-1< t \leq 0$ is  the coordinate in the direction of the inner normal. In particular, points of $\Gamma$ will have coordinates of the form  $(x',0)$ which we identify with $x'.$
At every point $x' \in \Gamma$ we split the cotangent space $T^*_{x'}(\R^n) \simeq \R^n$ into a direct sum $$T^*_{x'}\R^n= T^*_{x'}( \Gamma)\oplus \R\eta,$$ where $\eta$ is the conormal at $x',$ and use the coordinates on $T^*_{x'}\R^n$ of the form  $(\xi'_1,\dots,\xi'_{n-1},\nu),$ where $\xi' =(\xi'_1,\dots,\xi'_{n-1}) $ are the coordinates of a vector in $ T^*_{x'}(\Gamma)$ and $\nu$ is the coordinate along the conormal. 
 
Since  $ \cL_\la(x,D)$ is  properly elliptic,  $km=2r$ and,  for any $\la\in S^q,$   $x' \in\Gamma$ and any cotangent vector to the boundary $\xi' \neq 0 $  at $x',$  the polynomial  $P(\nu) =\mathrm{det}\, p (\la,x',\xi', \nu )$  has exactly $r$ roots in the upper half-plane $ \Im z > 0$  \cite{[Pe]}.

 In terms of the ordinary differential operator  $ p(\la, x',\xi', D_t ),$ obtained  by substituting  $\nu$ with   $D_t=i^{-1}{ \frac{d}{dt}}$ in the  principal symbol,   the above condition means that the  subspaces of stable (resp. unstable)  solutions  $ M^\pm (\la,x',\xi' )\subset  L^2( \mathbb R_\pm ;\C^m),$   whose elements are  solutions of the system $ p(\la,x', \xi', D_t)v(t)=0$ exponentially decaying to $0$ as $t \ra +\infty$  (resp. $t \ra -\infty$) are  both of dimension $r.$ In particular   $ M^\pm (\la,x',\xi' )$ are the fibers of  two  vector bundles  over $S^q\times [T^*(\Gamma) - \{0 \}],$ which will be denoted with $M^\pm.$ 

Let us denote with $\ga_j u$ the restriction  to  $\Gamma$ of  the function $ D^{(j)}_t u(x',t).$ 

Using the coordinates on $\cN, $ we rewrite the boundary operators  $\cB^i_\la $   in the form   \begin{equation}\label{bop}
 {\mathcal B}^i_\la(x,D)u=\sum_{j=0}^{k_i}  \cB^i_{j }(\la, x' ,D) \ga_ju,
\end{equation} 
where,  $\cB^i_{j }(\la, x',D)$ is a differential operator  of order $k_i -j$ acting on vector functions defined on the manifold $\gG .$  
In the new coordinates  the principal symbol of the boundary operator $\mathcal B_\la(x,D)$  
is  the matrix function $p_b(\la,x,\xi) $  whose $ i$-th row  is 
\begin{equation} \label {boundsymb}
p_b^{i }(\la,x,\xi)  =  \sum_{j=0}^{k_i}p_b^{ij}(\la,x',\xi') \nu^j,\end{equation}
where $p_b^{ij}(\la,x',\xi')$ is the principal symbol of  the operator $\cB^i_{j }(\la, x',D).$

By $H_1$ and $H_2,$ for any $\la \in S^q,$ the boundary  operator  $\cB_\la(x,D)$ verifies the  {\em Shapiro-Lopatinskij  condition} with respect to ${\cL}_\la( x,D).$ This means that,   for each $x' \in  \partial \Omega$   and each $\xi'$ belonging to $T_x'\partial \Omega,$ the subspace $ M^+ (\la, x',\xi' )$  is isomorphic to $\C^r$ via the map $b(\la, x', \xi')$  defined by 
\be \label{b1} b(\la, x', \xi') v = [ p_b(\la, x',\xi',D_t)v](0).\end{equation}  In particular  both  vector bundles $M^\pm$ are trivial.
\sk 
Identifying an endomorphism of $\C^r$ with its matrix in the  canonical basis we   define our second map 
$\tau \colon S^q\times [T^*(\Gamma) - \{0\}] \ra GL(r;\C)$  
by 
\be \label{partsig}  \tau(\la,x',\xi') = b(\la, x', \xi')b^{-1}(\infty, x', \xi'). \end{equation}

We will define degree of the matrix functions $\sigma$ and $\tau$  using Fedosov's approach in \cite{[Fe]}. For this we will need matrix-valued differential forms, or equivalently, matrices having differential forms as coefficients. The product of two matrices of this type is defined in the usual way, with the product of coefficients given by the wedge product of forms.   The matrix of differentials  $(d\sigma_{ij})$ will be denoted by $d\sigma.$ 

Let us consider  a compact oriented manifold $V$ of odd dimension $2v-1$   and a smooth map $\phi \colon V \ra GL(l;\C).$  Taking the trace of the  $(2v -1)$-th power of the matrix  $ \phi^{-1} d\phi $ we obtain an ordinary $(2v-1)$-form  $ tr (\phi^{-1} d\phi )^{2v-1} $ on $V.$  The {\it Bott-Fedosov degree}  of $\phi$ is defined by

   \begin{equation} \label{bottfed}
\deg(\phi)= \displaystyle N \int_{V} tr (\phi^{-1}d\phi)^{2v-1}, \end{equation}
where $N= - \displaystyle {\frac{(v-1)!} {(2\pi i)^{v}( 2v-1)!}}.$
   % As we will see later, Bott's integrality theorem implies that $\deg(\phi)\in \Z.$
    
\sk     
    
We define now  the interior and  boundary  multiplicity of the linearized  boundary value problem \eqref{lin}.

Let  $q$ be even, we associate to the $GL(m;\C)$-valued  function $\sigma$ constructed in \eqref{symbo}  the  one form $ \sigma^{-1} d \sigma$ defined on $S^q \times (\R^{2n} - K\times\{0\}).$   Without loss of generality we can assume that $\bar\Omega \times\{0\}$ is contained in the unit ball $B^{2n} \subset \R^{2n}$ and hence we can consider the restriction of  $ \sigma^{-1} d\sigma$  to  $S^q\times S^{2n-1} $ as a  one form on the compact manifold  $S^q\times S^{2n-1}.$ To be precise,  the latter is the pullback of the former by the inclusion   $i \colon S^q\times S^{2n-1} \ra S^q \times (\R^{2n} - K\times\{0\}).$  Being homogenous, $\sigma$ is uniquely determined by its restriction to   $S^q\times S^{2n-1}.$ Thus we will not distinguish in the notation the map $\sigma$ from its restriction  to  this space. On the other hand, the chain rule allows us  to denote  with  $ \sigma^{-1} d\sigma$  the pullback form too.  

By definition, the {\it interior  multiplicity}  of the family $(\cL,\cB)$ is
\begin{equation} \label{fed}
   \mu_i (\cL,\cB)= \deg(\sigma)=  \displaystyle {\frac{-(\frac{1}{2}q+n-1)!} {(2\pi i)^{(\frac{1}{2}q+n)}( q+2n-1)!}} \int_{S^q\times S^{2n-1}} tr (\sigma^{-1}d\sigma)^{q+2n-1}.
  \end{equation}

The   {\it boundary  multiplicity}   $ \mu_b (\cL,\cB)$   defined  in a similar way. Namely: 
   \begin{equation} \label{parfed}
 \mu_b (\cL,\cB)=  \deg(\tau)= \displaystyle {\frac{-(\frac{1}{2}q+n-2)!} {(2\pi i)^{(\frac{1}{2}q+n-1)}( q+2n-3)!}} \int_{S^q\times S(\gG)} tr (\tau^{-1}d\tau)^{q+2n-3},
  \end{equation}
 where  $S(\gG)$ is the unit sphere bundle of the cotangent bundle  $T^*(\gG).$
  
 It  follows from Fedosov's formula for the Chern character  of a family of elliptic pseudo-differential operators on $\R^n$ \cite[Corollary 6.5 ]{[Fe-1]}  and Bott's Integrality Theorem that   $\mu_i (\cL,\cB)$ is  an integer.   We will show in Section $4$ that $\mu_b (\cL,\cB)\in Z $ as well.
 
 Finally, the  {\it multiplicity} of $(\cL,\cB)$  is defined as 
\begin{equation}\label{multiplbif} \mu (\cL,\cB) = \mu_i (\cL,\cB)+\mu_b (\cL,\cB). \end{equation}

By construction, the integral number  $\mu (\cL,\cB)$  depends only on the principal symbols of the interior and boundary operators  and is invariant under homotopies of families of  elliptic  boundary value problems.

\begin{remark}  {\rm While the reality condition \eqref{real}  on the principal symbols of the interior and boundary operators may eventually place some restrictions on  the possible values of  the degree, the above condition was nowhere used  in the definition of  $\deg(\sigma)$ and $\deg(\tau).$  Hence exactly the same formulas allows to define the multiplicity  $\mu (\cL,\cB)$  of any family of linear elliptic boundary value problems with complex coefficients.}\end{remark}

Now, let us introduce the natural numbers $n(q).$  Denoting with $\nu_p(n)$ the power of the prime $p$ in the prime decomposition of $n\in \N,$ let $m$ be the  number-theoretic function constructed as follows: the value  $m(s)$ is defined  through its  prime decomposition,  by setting for  $p=2,$   $ \nu_2\left(m(s)\right)=2+\nu_2(s)$ if  $s\equiv 0 \mod 2$ and $ \nu_2\left(m(s)\right)=1$ if the opposite is true. While, if $p$ is an odd prime, then
  $\nu_p\left(m(s)\right)= 1 +\nu_p(s)$ if  $s\equiv 0 \mod (p-1)$  and $0$ in the remaining cases. In particular  $m(s)$ is always even.  The function $m$ was introduced by Adams \cite{[Ad]}.  It is well known that for  $q=4s$  the group $ J(S^q)$ is a cyclic group of order $m(2s).$  
  
 For $q\equiv 0,4 \mod 8,$   let  \begin{equation} \label{muca}  n(q)= \begin{cases} m(q/2) & \hif \, q\equiv 0 \mod 8\\  2 m(q/2) & \hif \, q\equiv 4 \mod 8 \end{cases} \end{equation} 
and $m$ is the number theoretic function defined above.

With all of the above said we can  state our criteria for bifurcation of solutions of  \eqref{bvp1}. 
  
 \begin{theorem}\label{th:40}
Let the problem \begin{equation}\label{bvp2}
\left\{\begin{array}{l} \cf\,(\la, x,u,\ldots ,D^{k}u) = 0,\, x \in  \Omega \\  \cg(\la, x,u,\ldots, D^{k_i} u)=0,\, x \in  \partial \Omega , \end{array}\right. \end{equation}  verify   assumptions $ H_1, H_2$ and $H_3.$ 
 If $q\equiv 0,4 \mod 8,$ then  bifurcation of smooth solutions of \eqref{bvp2}  from some  point  of the  trivial branch  arises provided that 
 $\mu( L,B)$
 is not divisible by $n(q).$
 \end{theorem}
  
  \begin{remark}{\rm If the principal part of the boundary operator is independent of $\la,$  then $\mu_b(L,B)=0$ and we obtain Theorem 1.4.1  of \cite{[Pe]}. If instead the principal part of the operator  $\cL_\la(,x,D)$ is independent of $\la,$ then bifurcation of solutions is determined  by  $\mu_b(L,B),$ i.e., by Bott-Fedosov  degree of  $\tau$.}  \end{remark}

If $(\cf,\cg)$   verifies   assumptions $ H_1, H_2$ and $H_3$ of Theorem \ref{th:40}, then for  any lower order perturbation  \[( \cf'(\la, x,u,\ldots ,D^{k-1}u), {\cg'}^i(\la, x,u,\ldots, D^{k_i-1} u)),\, 1\leq i \leq r, \] such that  $\cf'(\la, x,0)=0,$ ${\cg'}^i(\la, x,0)=0 $ and  such that the coefficients of the linearization of $(\cf',\cg')$ converge  uniformly  to $0$ as $\la\ra \infty,$ also  the perturbed  problem 
\begin{equation}\label{bvp4} \left\{\begin{array}{l} \cf(\la, x,u,\ldots ,D^{k}u) +  \cf'(\la, x,u,\ldots ,D^{k-1}u) = 0,\, x \in  \Omega \\  \cg^i(\la, x,u,\ldots, D^{k_i} u) + {\cg'}^i(\la, x,u,\ldots, D^{k_i-1} u) =0,\, x \in  \partial \Omega, \end{array}\right. \end{equation}
 verifies the same assumptions. Taking into account that lover order perturbations do not affect the value of $\mu(\cL,\cB)$ we have: 
 
 \begin{corollary}\label{lowerorder}  If $(\cf,\cg)$ verifies all of the  assumptions in Theorem \ref{th:40} and if the perturbation $(\cf',\cg')$ is as above, then 
 there must be some bifurcation point $\la \in \R^q$   for solutions of  the perturbed problem \eqref{bvp4}  as well.\end{corollary}
 
  The above corollary  shows that the bifurcation criterium based on the invariant $ \mu(\cL,\cB)$ is  robust.  Bifurcation invariants of local type \cite{[Iz]} lack of  this property.  On the other hand   Theorem  \ref{th:40}  does not give any information about where  the bifurcation points are located. 
  
  In many instances   $ \mu(L,B)$  vanishes, e.g., when the top order terms of 
  the operator  $(\cLL ,\cB_\la )$ is independent of  $\la.$ However, examples of boundary value problems with non vanishing $\mu(L,B)$ will  be constructed  in Section $5$ taking into account  the topology of $GL(m; \C).$ In some simple cases  the bifurcation index  can be computed from the coefficients of the linearization, using  rather  elementary index theorems, see \cite{[Pe-1]}  and \cite{[Pe-2]}.   
  \sk

\section{The index of bifurcation points } \label{sec:2}
Here we will shortly recall the definition of the  index of  bifurcation. The construction in \cite{[Pe]} uses  the index bundle  of a family of Fredholm operators  and  the generalized $J$-homomorphism.

If $\Lambda$ is a compact space, the set $Vect(\Lambda )$ of all isomorphism classes of real vector bundles over  $\Lambda$ is a semigroup under the direct sum. The  (real) Grothendieck group $KO(\Lambda )$ of a compact topological space  $ \Lambda$    is the quotient of the  semigroup $ Vect(\Lambda )\times Vect(\Lambda )$ by the diagonal sub-semigroup. Elements of $KO(\Lambda )$  are called  virtual bundles.  Each virtual bundle can be written  as a difference  $[E] - [F],$ where $E, F$ are vector bundles over  $\Lambda  $ and $[E]$ denotes the equivalence class of $(E,0).$  It is easy to see  that $ [E] - [F]= 0$ in $KO(\Lambda )$ if and   only if $E$ and $F$  become isomorphic after the addition of a trivial vector bundle  to both.  The complex Grothendieck  group $K(\Lambda)$ is defined by taking complex vector bundles instead of the real ones.  In what follows the trivial bundle with fiber  $\Lambda\times V$ will be denoted by  $\Theta(V).$ The trivial bundle  $\Theta(\R^n)$ will be simplified to  $\Theta^n.$ 
 
    Let $X,\ Y$ be real Banach spaces and let  $L\colon  \Lambda\rightarrow \Phi(X,Y)$     be a continuous family of Fredholm operators. Since $\coker L_\la$ is finite dimensional, using  compactness of $ \Lambda,$ one can find a finite dimensional subspace  $V \hof\,  Y$ such  that
             \begin{equation} \label{1.1}
\hbox{\rm Im}\,L_\la+ V=Y  \ \hbox{\rm for all }\  \la \in \Lambda.
\end{equation}

Because of the  transversality condition \eqref{1.1}  the family of finite dimensional subspaces 
$E_{\la}=L_{\la}^{-1}(V)$ defines a vector bundle $E$ over $ \Lambda.$ 

The {\it index bundle} of the family $L$  is the virtual bundle
 \begin{equation} \label{defind}  \Ind L= [E]-[\Theta (V)] \in KO(\Lambda ).\end{equation}
  
The index bundle has the same properties of the numerical index. It is homotopy invariant and hence invariant under perturbations by families of compact operators. It is additive under direct sums and  the same holds for the index bundle of the family of composed operators (logarithmic property). Clearly $\Ind L=0$ if $L$ is homotopic to a family of invertible operators. Below we will use the above  properties without any further reference. 
 
 Notice that the index bundle of a family  of Fredholm operators of index $0,$ belongs to the reduced  Grothendieck group $ \tilde{KO}(\Lambda )$  defined as  the kernel of the rank homomorphism  $ \rk \colon KO(-) \ra \Z,$ $rk([E]-[F]) = \rk E -\rk F.$

 Given a  vector bundle $E,$ let  $S[E]$ be the associated unit  sphere bundle with respect to some  chosen scalar product on $E$. Two vector bundles $E,F$ are said to be {\it stably fiberwise homotopy equivalent} if, for some $n,m,$ the unit sphere bundle $S ( E\oplus \Theta^n)$ is fiberwise homotopy equivalent to the unit sphere bundle $S (F\oplus \Theta^{m}).$ Let   $T(\Lambda )$ be the subgroup of $\tilde{KO}(\Lambda )$  generated by elements  $ [E]-[F] $ such that $E$ and $F$ are  stably fiberwise homotopy equivalent.  Put  $\J(\Lambda )= \tilde {KO}(\Lambda )/ T(\Lambda ).$ The projection to the quotient $\J \colon\tilde {KO}(\Lambda )\rightarrow \J(\Lambda )$ is   {\it  the generalized $\J$-homomorphism}. 
 
 \begin{remark} { \rm The group $\J(\Lambda )$  was  introduced by Atiyah in \cite{[At-1]}. He proved  that  $\J(\Lambda )$  is a finite group if $\Lambda$ is a finite $CW$-complex by showing  that  $\J(S^n)$ coincides with the image  of the stable  $j$-homomorphism of G. Whitehead.} \end{remark}

Now, let us introduce the index of bifurcation points constructed in  \cite{[Pe]} for    families of  Fredholm maps of index $0$  having as  range a  Kuiper space $Y$, i.e.,  a Banach space whose group of linear automorphisms is contractible.
  
  Let $U$ be an open subset of a finite CW-complex $ \Lambda$ and  let $O$ be an open subset of a Banach space X. Let $f:U\times O \rightarrow Y$ be a family of $C^1$ Fredholm maps of index 0  such that $f(\la ,0)=0$ \cite{[Pe]}.We will denote by $ L$  the family  $ \{Df_\la (0); \la \in U \}.$   
  
  The  pair $(f,U)$  is called  {\it admissible}  if  the  set $\Sigma (L)= \{ \la\,  |\, \Ker L_\la \neq 0\}$   is a  compact proper subset of $U.$
 
  If $(f,U)$ is admissible, the {\it  index of bifurcation points of $f$ in $U$ } is defined by 
 \begin{equation}
\label{eq:2.6}
\beta (f,U)=\J(\Ind \bar L),
\end{equation}
where $\bar L\colon \Lambda \rightarrow  \Phi_0(X,Y)$  is any family which coincides with $L$ 
on an  open  neighborhood of $\Sigma (L)$ with compact closure contained in $U.$

It was shown in \cipe that $ \beta (f,U)$ verifies  the  homotopy invariance, additivity and excision property in their usual form. Here we are mainly  interested in: 
\vskip5pt

{ \it Existence  property:} \  If $\beta(f,U)\neq 0,$ then the family $f$ has  at least one   bifurcation point $\la_*$ in $U,$ i.e., a point such that  every neighborhood $V$ of  $(\la_*,0)$ contains a solution of $f(\la,x)=0$ with $x\neq 0.$
\vskip5pt

{\it Normalization property:} \  $\beta (f,\Lambda)=\beta(f) =\J(\Ind L)$. 
\vskip5pt

\sk

\section {Proof of  theorem \ref{th:40} .} 

   Denoting with $\hbbr$ the product $ \prod_{i=1}^r H^{k+s-k_i-1/2}(\partial \Omega; \R),$
 the family of  nonlinear differential operators  \[(\cf,\cg) \colon \R^q \times C^\infty( \bar\Omega;\R^m) \rightarrow C^\infty( \bar\Omega;\R^m) \times\ C^\infty(\partial\Omega;\R^r)\]
  induces a smooth map 
\begin{equation}
\label{nbvp1}
 h=\colon \R^q \times H^{k+s}(\Omega;\R^m ) \rightarrow    \hsbr
\end{equation}
having   $\R^{q} \times\{0\}$ as a trivial branch (see \cite[Section 5.2]{[Pe]}). 

The Frechet derivative $Dh_\la(0)$ of the map  $h_\la$ at  $u\equiv 0$  coincides with the  operator 
$(L_\la,B_\la)$ induced by the linearization $(\cL_{\la}, \cB_\la).$  Since, for any $\la\in R^q,$ $(\cL_{\la}, \cB_\la)$ is elliptic, using proposition 5.2.1 of  \cite{[Pe]}, we can find a neighborhood  $O$  of $0$ in $\hkr$  such that $h\colon \R^{q} \times O \rightarrow  \hsbr $  is a smooth  family of  semi-Fredholm maps. 

 By  $H_2,$  the family  of boundary value problems $(\cL, \cB)$ extends to a smooth family  parametrized by $S^{q},$ which will be denoted in the same way. Let $(L,B)$ be  the family of operators induced on Hardy-Sobolev spaces. 
   
  Being   $(L_{\infty},B_{\infty})$   invertible,  by continuity of the index of semi-Fredholm operators, $(L_{\la},B_{\la})$  is Fredholm of index $0$ 
   for all $\la\in S^{q}$ which, on its turn, implies that  the map  $h\colon \R^{q} \times O \rightarrow  \hsbr $  is a smoothly parametrized family of Fredholm maps of index $0.$  
   
     Since  $(L_\la,B_\la)$ is invertible in a neighborhood of  $\infty\in S^q,$ $\Sigma (L,B) $ is a compact subset of $\R^q\subset S^q.$  Therefore,  the pair $(h,\R^q)$ is admissible and  the index of  bifurcation points $\beta(h,\R^q)$ is defined. On the other hand,  $(L,B)$  is  an extension of the family $\la \ra  Dh_\la (0)$ to all of $S^q.$ By  definition of the index of bifurcation, $\beta(h,\R^q)= \J(\Ind(L,B)).$  

If, under the hypothesis of Theorem  \ref {th:40}, 
 we can show that $\J(\Ind(L,B)) \neq 0 $ in $\J(S^q),$  then  the family $h$ must have   a bifurcation point  $\la \in \R^q,$ by the existence property of the bifurcation index.  This would complete the proof of the theorem,  since by Proposition 5.2.2 of \cipe  a bifurcation point of the map  $h$  is also a bifurcation point for smooth classical solutions of  \eqref{bvp1} in the sense of definition \ref{def:1}.

We will show that  $\J(\Ind(L,B)) \neq 0 $ in $\J(S^q)$ in three steps. 

\vskip 10pt
{\bf Step 1} We  assume first  that the principal part of $\cB_\la$ is independent of $\la$ 
and hence $\mu(L,B) =\deg(\si).$
 This is precisely the case considered in \cite[Theorem 1.4.1]{[Pe]}. In that paper, we showed that, if $\deg(\si)$ is not divisible by $n(q),$ then $\J(\Ind L) \neq 0.$

\sk

{\bf Step 2}  Let us  assume now that the principal part of  $\cLL$ is independent of $\la$ while there are no restrictions on the boundary  condition $\cB.$  

 In this case we will use a family version of the Agranovich-Dynin reduction (see \cite{[Ag-Dy], [W-R-L]}). In the case without parameters this reduction  computes the difference between the indices of two elliptic boundary value problems  for the same interior  operator $\cL(x,D)$ as the index of a pseudo-differential operator on the boundary $\Gamma$ whose principal  symbol  is constructed  in terms of the data.  The  result extends easily to families of elliptic boundary value problems.  While discussing  the Agranovich reduction  in \cipe  we  provided full details, we will be slightly  more sketchy here.

We will consider only  classical   pseudo-differential operators acting  on complex vector functions on a compact smooth orientable manifold $\Gamma.$ This is exactly the same class of pseudo-differential operators as the one introduced in \cite{[At-Si-I]} but we are dealing here with trivial bundles only. For  a detailed exposition see  \cite{[R-S],[W-R-L]}. 
The class of  pseudo-differential operators  under consideration  contains all differential operators on manifolds  and  moreover is invariant under composition and formation of adjoints. We will need only a few facts about  this class.
\sk 
1) For  $s\in \N$, the Hardy Sobolev space $\hsman{\Gamma}$ can be defined as  the space  of all vector functions such that $\cL(x,D) u \in L^2(\Gamma;\C^m)$ for every differential operator of order less or equal than $s.$ This definition extends to all $s\in \R$ in the usual form, using the square root of the Laplacian. Every  pseudo-differential operator $\cR$ of order $k$ defines a bounded operator  $R\colon \hkman{\Gamma} \rightarrow \hsman{\Gamma}.$

2) Each  pseudo-differential operator $\cR$ of order $k$ on $\Gamma$ has a well defined principal symbol $ \rho_k\colon T^*(\Gamma)-\bar 0 \ra \C^{m\times m},$ where $\bar 0$ denotes the image of the zero section. The principal symbol $\rho_k$ is a
 smooth, positively  homogeneous function of order $k.$ It  behaves well under composition of operators. 
Namely, $ \rho_{k+s}(\cR\cQ)= \rho_k(\cR) \rho_s(\cQ).$

3) A pseudo-differential operator $\cR$  of order $k$  is  elliptic  if  $\rho_k(\cR)(x,\xi) \in GL(m;\C)$ for all $(x,\xi)$ with $\xi \neq 0.$  Every elliptic operator  possesses a (rough) {\it parametrix},  which is a proper  pseudo-differential operator  $\cP$   of order  $-k$ such that  both  $\cR \circ\cP -\Id $ and  $\cP\circ\cR -\Id$ are of order $-1.$
Since operators of order $-1$ viewed as operators of $\hsman{\Gamma}$ into itself are compact, by the Riesz characterization, the bounded operator induced by an elliptic pseudo-differential operator on Hardy-Sobolev spaces is Fredholm. 
\sk
 
The Agranovich-Dynin reduction  can be carried out in the context of continuous families of elliptic boundary value problems.  Although here we will consider smooth families only,  since we have defined  the
bifurcation index  for continuous families of $C^1$-Fredholm maps we will  work out the reduction in the above generality.   Continuous families of pseudo-differential operators have been introduced 
in \cite{[At-Si-IV]} (see also  \cipe for operators on $\R^n$). 

Let us denote by $S(\Gamma) \subset T^*(\Gamma) $ the unit sphere bundle with respect to the induced riemannian metric on $\Gamma.$    We will need also  the following:

\bl\label{app} 
 Every continuous map $\rho \colon \Lambda \times S(\Gamma) \ra Gl(m,\C)$ can be uniformly approximated by  the restrictions to  $\Lambda \times S(\Gamma)$ the principal symbols  $\rho(\cS)$  of continuous families of pseudo-differential operators $\cS$ of any chosen order. 
 Moreover, if $\Lambda$ is a smooth manifold   the approximating symbols can be chosen  smooth (in the parameter variables as well). \el
\proof   The first assertion is proved in \cite[Proposition 6.1]{[At-Si-IV]}. The second  is an immediate consequence of the method in proof there. 

  \vskip10pt

Let us  consider  two families  of linear elliptic  boundary value problems $(\cL,\cB_+)$ and $(\cL,\cB_-)$ with the same constant 
 interior operator $\cL ,$ and with parameter dependence  only on the boundary operators: 

\be\label{ellitt}
\left \{\ba{lll}\mathcal L(x,D)&=&\sum_{|\alpha|\leq k}a_{\alpha }(x)D^{\alpha },\\ 
\mathcal{B_\pm}^i(\la, x,D)&=&\ga_0\sum_{|\alpha|\leq k_i}{b^\pm}^i_\alpha(\la,x)D^{\alpha},\, 
1\leq i\leq r. 
\ea \right.
\end{equation}
Here the  parameter $\la$ belongs  to a compact connected topological space $\Lambda,$
while the  matrix functions  $a_{\alpha }(x) \in \C^{m\times m}, \, {b^\pm}^i_\alpha(\la,x)  \in \C^{1\times m}$  are smooth in $x,$ and depend continuously on $(\la,x)$ together with all their partial derivatives.

Using the coordinates $(x',t)$ on $\cN, $ we rewrite the boundary operators  $\cB_\pm^i(\la, x,D)$  in the form  described in  \eqref{bop}, i.e., 
\begin{equation}\label{bop'}
 {\mathcal B_\pm}^i(\la,x,D)u=\sum_{j=0}^{k_i}  {\cB_\pm}^i_{j }(\la, x' ,D)\ga_ju.
\end{equation} 

Here,  however, we will have to extend our  considerations  to  pseudo-differential boundary conditions. Namely,  the coefficients ${\cB_\pm}^i_{j }(\la, x',D)$ can be  pseudo-differential operators on $\Gamma$  of order $k_i -j.$ The principal symbol of the boundary operator is defined in the same way as  in \eqref{boundsymb} and there are no changes in the formulation of the Shapiro-Lopatinskij  condition. 

\sk

 For $(\la,x',\xi') \in \Lambda \times (T^*(\Gamma)-\bar 0),$ let us define
\be \label{partsig2}  \tau(\la,x',\xi') = b_+(\la, x', \xi')b_-^{-1}(\la, x', \xi'), \end{equation}
where $$b_\pm(\la, x', \xi')\colon  \cM^+ (\la, x',\xi' )\ra \C^r $$ is the isomorphism  associated by  \eqref{b1} to the boundary operator ${\cB_\pm}(\la, x,D).$

  The following proposition is the Agranovich-Dynin reduction  for families:
   
\begin{proposition}\label{adred} 

Given $(\cL, \cB_\pm)$ as above, there exists a family $\cS$   of pseudo-differential operators  of order $0$ on $\Gamma$ such that 

         \be\label{ad} \Ind (L,B_+) -\Ind(L,B_-) = \Ind S.\end{equation}
            
                Moreover the restriction of  the principal symbol $\rho_0(\cS)$ to $S^q\times S(\Gamma)$  can be taken arbitrarily close in the sup  norm  to the restriction of $\tau$ to the same subspace. 
       \end{proposition}
       
\proof 

 Using  lemma \ref{app} we can find a smooth  family $\cS = \{\cS_\la; \, \la\in  \Lambda\}$   of elliptic  pseudo-differential operators of order $0$ on $\Gamma$   such that the family of principal symbols restricted to $S^q\times S(\Gamma)$ is arbitrarily close to the restriction of $\tau.$    Since the Shapiro-Lopatinskij  condition  \eqref{b1}   is stable under small perturbations, it follows that  
the family $(\cL, \cS\cB_-)$ is a family of elliptic boundary value problems with pseudo differential boundary conditions.

Since  $(L,SB_-)=(\Id \times S)(L,B_-),$  by the logarithmic property of the index bundle,
   \[ \Ind(L,S B_-)= \Ind (\Id \times S) +\Ind (L,B_-)= \Ind S +\Ind (L,B_-)\]
   
  On the other hand, if on  $S^q\times S(\Gamma),$ $\rho_0(\cS)$ is close enough to $\tau,$  the affine homotopy 
  \[(1-t)(L, SB_-) +t(L, B_+)\] is a homotopy of  families of linear Fredholm operators between 
  $(L,SB_-)$   and $(L,B_+).$   The proposition now follows from the homotopy invariance property of the index bundle. \qed

\sk
\sk\sk

 \begin{theorem}\label{th:41}
Let the problem \begin{equation}\label{bvp5}
\left\{\begin{array}{l} \cf\,(\la, x,u,\ldots ,D^{k}u) = 0,\, x \in  \Omega \\  \cg_i(\la, x,u,\ldots, D^{k_i} u)=0,\, x \in  \partial \Omega , 1 \le  i \le  r ,\end{array}\right. \end{equation}  verify   the assumptions $ H_1, H_2$ and $H_3.$ Assume 
moreover that the family of  interior operators $\cL(x,D)$ of  the linearization at $u\equiv 0$ is independent of $\la.$ 

If $q\equiv 0,4 \mod 8,$  there exists  at least one bifurcation point from the trivial branch   provided that $\mu_b(\cL,\cB) $ is not divisible by $n(q).$                                                                            
\end{theorem}

\proof 

We are going to compute $\J (\Ind (L,B))$ from  $\mu_b(\cL,\cB) =\deg(\tau)$ using  the complexification $(L^c,B^c)$ of the linearized equations at $u\equiv 0.$

Since  $\ker( L^c,B^c) =\ker (L,B)\otimes \C ,$  from  the
definition of the index bundle in \eqref{defind}  it follows that 
\be\label{compl} \Ind (L^c,B^c) = c(\Ind(L,B)),\end{equation}  where  $c\colon \tilde{KO} \rightarrow \tilde K$ is the complexification homomorphism.
 
By Bott periodicity, for $q=4s,$ both $\tilde K(S^q) \cong \Z$  and  $\tilde{KO}(S^q) \cong  \Z$ are infinite cyclic with  
$\tilde K(S^q) $ generated by  powers   $\xi_q =\left ([H] -[\Theta^1]\right)^{2s},$ where $H$ is the tautological line  bundle over $P^1(\C).$ 
 Moreover, by  \cite[section 13.94]{[Sw]}  $c\colon\tilde{KO}(S^q) \rightarrow \tilde{K}(S^q)$ is an isomorphism for $q\equiv 0 \mod 8 $  and a monomorphism with image generated by $2\xi_q$ for  $q\equiv 4 \mod 8.$ 
 
   We choose as generator of $\tilde{KO}(S^q)$ an element $\nu_q$ such that \begin{equation}\label{gen} c(\nu_q)= \begin{cases} \xi_q  &\hif q\equiv 0 \mod 8 \\ 2\xi_q  &\hif  q\equiv 4 \mod 8.\end{cases}\end{equation} 
  Then each  element $\eta \in\tilde K(S^q)$ with $q=4s$ is  uniquely determined by its  {\it degree}  $d(\eta) \in \Z$  verifying  $ \eta = d(\eta)\, \xi_q, $ and
   each  element  $ \eta$ of $\tilde{KO}(S^q)$  has  a degree defined in the same way.  Clearly,  
\begin{equation}\label{degr} d(c(\eta)) = \begin{cases} d(\eta) &\hif q\equiv 0 \mod 8 \\ 2d(\eta)  &\hif  q\equiv 4 \mod 8.\end{cases}\end{equation}
 
Let us denote with  $\rm{H}^{*}(-;\C)$  the de  Rham cohomology with coefficients in $\C$ and compact supports.  We will denote with  $\rm{H}^{ev/odd}(-;\C)$ the cohomology  in even degrees and odd degrees respectively.  By the uniqueness of the Chern character and Bott's integrality theorem (\cite[Theorem 9.6, Chap.18]{[Hu]}),  $\ch \tilde{K}(S^q) \ra \rm{H}^{ev}(S^q;\C)$ sends  $\tilde{K}(S^q)$  isomorphically into $\rm{H}^{ev}(S^q;\Z) \subset \rm{H}^{ev}(S^q;\C).$  Hence  the degree of an element  $\eta \in \tilde{K}(S^q)$ can be computed by  evaluating the  Chern character %$ \ch \colon \tilde{K}(S^{4s}) \rightarrow  \rm{H}^{ev}(S^{4s};\C)$
 on the fundamental class $[S^{q}]$   of  the sphere. Namely, 
\begin{equation}\label{cher} d(\eta) = < \ch( \eta); [S^{q}]> . \end{equation}  

 We  will compute the degree  of   $ \Ind (\cL^c,\cB^c )$ using  \eqref{cher}  and  the Agranovich-Dynin  reduction.

Put  $(\cL, \cB_+) =  (\cL^c, \cB^c)$ and $(\cL, \cB_-) =  (\cL^c, \cB^c_\infty)$
in  Proposition \ref{adred}.  Being  $(\cL^c, \cB^c_\infty)$ a constant family,  $\Ind (\cL^c, \cB^c_\infty)=0,$ and hence by Proposition \ref{adred}  we have,
 \begin{equation} \label{indsing}  \Ind (L^c,B^c) = \Ind S, \end{equation}
where $S$ is induced  by a family of pseudo-differential operators $\cS$ on $\Gamma $ whose principal  symbol $\rho  =\rho_0(\cS)$  is homotopic to $\tau.$
\sk
   Now let us  apply the cohomological form of the Atiyah-Singer family index theorem to $\cS.$  
   
   Since $\Gamma $ is a boundary, its Todd class vanishes. Thus by the Atiyah-Singer family index theorem \cite[Theorem (3.1)] {[At-Si-IV]}  
   \begin{equation}\label{as} \ch (\Ind S) =  p_* \ch[\rho],\end{equation}   where $p_*\colon H^*( S^q\times T^*(\gG) ) \ra  H^*(S^q) $ is the direct image homomorphism associated to the bundle of tangents along the fiber and  $[\rho] \in \K_c( S^{q}\times T^*(\Gamma))$  is  the symbol class of $\rho.$ 
 
By definition, the {\it symbol class} $[\rho] \in \K_c( S^{q}\times T^*(\Gamma))$ is  obtained  from the map $\rho\colon S^q \times (T^*(\Gamma)-\bar 0)\ra GL(r;\C) $ by means of the  {\it clutching  construction} \cite{[At-Si-I]} described below: 

Let $S^*(\Gamma)$ be the  fiberwise compactification  of $T^*(\Gamma),$ obtained by adjoining a  point at infinity to each fiber, and let  $S^*= S^q \times S^*(\Gamma).$
 Then $S^*$ is the union of two open sets $U_0=S^q \times T^*(\Gamma)=S^q \times (S^*(\Gamma)-\bar\infty )$ and 
 $U_1 = S^q\times (S^*(\Gamma) - \bar 0),$ where as before  $\bar \infty $ denotes  the image of the section at infinity. Let $D_i  \subset U_i$ be the set of points $(\la,v)\in S^*$ with the  norm of $||v|| \leq 1$ and $||v|| \geq 1$ respectively.  
 
    We obtain a  vector bundle $E$ over $S^*$   gluing two trivial bundles $\theta^r$ with fiber $\C^r$  over $D_i ,  \  i=0,1$ by means of the restriction of $\rho$ to $ S^q\times S(\Gamma).$  Since the restriction of $E$  to a neighborhood of $\bar\infty $ is trivial,  $ [E]-[\Theta^r]$  defines an element $ [\rho]$ belonging to  $\tilde K(S^*/ \bar\infty)\cong  \K_c( S^q\times T^*(\Gamma)).$ By definition,  the above element is the {\it symbol class} of $\rho$. 
    
 The symbol class $ [\rho]$
 is defined  in terms of the restriction of $\rho $ to $S^q\times S(\Gamma)$ only.  Indeed,  the above construction associates a  (homotopy invariant) symbol class $ [\rho]\in   \K_c( S^{q}\times T^*(\Gamma)) $  to any continuous map   $\rho \colon S^q\times S(\Gamma) \ra GL(r;\C).$ 
    
   \begin{remark}  {\rm The formula  \eqref{as} differs from the one in \cite[Theorem (3.1)] {[At-Si-IV]}  by a factor $(-1)^{n-1}.$  This factor, which is irrelevant to our considerations, disappears  by substituting  the  orientation of $ T^*(\gG)\otimes \C$ used  in the above paper with  the one  in \cite[ Theorem 2, Chap XIX]{[Pa]}.}\end{remark}

  \sk 
  In \cite [¤3]{[Fe]},   Fedosov chooses  two trivializations  of $E_{| U_i}$ whose   transition function over $U_1\cap U_2$  coincides with  $\rho.$ Using this trivializations he defines a connection on $E$ and uses its  curvature  in order to construct  a  (non homogeneous) differential form representing the Chern character of $[\rho].$   The  result in  \cite [¤3, (17) ]{[Fe]} is  that  $\ch( [\rho] )$
 is the cohomology class of the differential form: 

\begin{equation} \label{cero}-\sum_{j=1}^\infty \displaystyle { \frac{(j-1)!} {(2\pi i)^j(2j-1)!}}  d\big(h(\|v\|) tr(\rho^{-1}d\rho)^{2j-1}\big),\end{equation}
  where $ h(t)$  is a smooth function which vanishes in a neighborhood of  $0$ and such that  $h(t)=1,$ for $t\geq1.$ 
   Actually, in \cite[¤3] {[Fe]} only the case $\Lambda=pt$ is considered, but all his  arguments hold word for word   for families parametrized by compact orientable manifolds.
   
On the other hand, since the restrictions  of $\rho$ and $\tau$ to $S^{q}\times S(\Gamma)$  are homotopic,  the vector bundles obtained by gluing trivial bundles using either  $\rho$ or  $\tau$ are isomorphic and therefore  their Chern characters coincide.  In conclusion we obtain 
\be\label{cernia} \ch( [\rho] )= \ch( [\tau] )=\, \big\{-\sum_{j=1}^\infty \displaystyle{ \frac{(j-1)!} {(2\pi i)^j(2j-1)!}}  d\big(h(\|v\|) tr(\tau^{-1}d\tau) ^{2j-1}\big)\big\},\end{equation}
where $\big\{ \theta \big\}$ denotes the cohomology class of the form $\theta.$ 
\sk

The direct image homomorphism $p_*$ in de Rham cohomology is  the homomorphism induced by a cochain  homomorphism  called {\it integration along the fiber}. The latter 
  takes $(d+2n-2)$-forms with compact support  on $S^{q}\times T^*(\Gamma)$ into $d$-forms on $S^{q}$ literally by integrating  the fiber variables. (see  \cite[ ¤ VII]{[G-H-V]}). 
  
 Denoting with  $\oint$ the integration along the fiber,  
from \eqref{cernia} we get %obtain that  $ p_*\ch ([\rho])$  is the cohomology  class of 

\be \label{cherone}  p_*\ch ([\tau])=\big\{  \sum_{j=n-1}^\infty \displaystyle { \frac{(j-1)!} {(2\pi i)^j(2j-1)!}} \oint_{T^*(\Gamma)}  d\left(h(\|v\|) tr(\tau^{-1}d\tau)^{2j-1}\right)\big\}. \end{equation}

 On the other hand,  the  evaluation of a cohomology class on the fundamental class of an $n$-manifold   in de Rham cohomology corresponds, at the cochain level, to   the integration of a representing form over the manifold.  Therefore, integrating over $S^q=S^{4s}$   the $4s$-homogenous term  from 
  \eqref{cherone} and using Fubini's theorem for integration along  the fiber  \cite[ ¤ VII]{[G-H-V]},  we get
\begin{equation}\label{cher1}   <p_*\ch ([\tau]); [S^{q}]>=   N \int_{S^{4s} \times T^*(\Gamma)}  d[h(\|v\|) tr(\tau^{-1}d\tau)^{4s+2n-3}],\end{equation} 
  where   the right hand side is the ordinary integration of  the $(4s+ 2n-2)$-form    over a manifold of the same dimension and 
   \begin{equation}\label{number} N=- \displaystyle { \frac{(2s+n-2)!} {(2\pi i)^{2s+n -1}(4s+2n-3)!}}.\end{equation} 
   It is easy to see that  $$ d[ tr(\tau^{-1}d\tau)^{4s+2n-3}] =-tr d[ (\tau^{-1}d\tau)^{4s+2n-2}] =0,$$ 
   and since $h(\|v\|) \equiv 1$  if $\|v\|\geq 1,$ the differential form $d[h(\|v\|) tr(\tau^{-1}d\tau)^{4s+2n-3}]$ vanishes outside $D_0. $ 
   
    Thus  \eqref{cher1} reduces to an integral over the manifold with boundary $D_0.$  Using Stokes theorem  we obtain 
  
  \begin{equation}\label{cher2}  <p_*\ch ([\tau]); [S^{q}]>=  N \int_{S^{4s} \times S^*(\Gamma)}  tr(\tau^{-1}d\tau)^{4s+2n-3}= \deg(\tau) = \mu_b(\cL,\cB)\end{equation} 

From \eqref{cher2},  \eqref{indsing}, \eqref{as} we have
  \begin{equation}\label{cher0}     < \ch \Ind(L^c,B^c) ; [S^{q}]> = < \ch \Ind S ; [S^{q}]> =
  \mu_b(\cL,\cB) .\end{equation}
 
 In particular,  by \eqref{cher0} and  Bott's Integrality  Theorem,  $\mu_b(\cL,\cB) \in \Z.$  
  
 Since  complexification of the index bundle of $\Ind(L,B)$  is the index bundle of the family of complexified operators  another consequence of \eqref{cher0} together with \eqref{cher}  is that  $d( c\big(\Ind (L,B)\big)=\mu_b(\cL,\cB).$ 
 
 From the above observation, using  \eqr{degr} we finally obtain 
\begin{equation}\label{degre} d\left(\Ind (L,B)\right) = \begin{cases}\mu_b(\cL,\cB) &\hif q\equiv 0 \mod 8 \\ \frac{1}{2} \mu_b(\cL,\cB)  &\hif  q\equiv 4 \mod 8.\end{cases}\end{equation}

 On the other hand, for  $q =4s,$  $J(S^{4s})\simeq  Z_{m(2s)}$ and  $J(\Ind L) =0$ if and only if  $d\left(\Ind( L,B)\right)$ is divisible by $m(2s).$ Hence  Theorem \ref{th:41} follows from \eqref{degre} and the definition of  $n(q)$ in \eqref {muca}.\qed 
 
 \sk
{\bf Step 3} We will reduce the general case to the two considered previously.  

Together with the  family of linearizations along the trivial branch  $(\cL,\cB)$ we consider $(\cL_\infty, \cB_\infty)$  as a constant family and compare  the following two  families   of elliptic  boundary value problems: 

\be\label{comp1}\left(\cL^1, \cB^1 \right)=\left(\cL_\infty  \cL ,\cB_\infty  \cL, \cB \right)\ \
 \end{equation} 
 and 
\be\label{comp2}\left(\cL^2, \cB^2 \right) =\left(\cL \cL_\infty, \cB_\infty \cL_\infty, \cB \right). \end{equation} 

\sk

Put  $X=H^{2k+s}(\Omega; \R^m ),\, Y=H^{k+s}(\Omega; \R^m ), \, Z=H^{(k+s)*}(\partial \Omega; \R^r) ,$  $V=H^{s}(\Omega; \R^m )$ and  $W =\hbbr.$
 
The operator 
$(L^1,B^1)\colon X \rightarrow  V \times W \times Z $
induced by \eqref{comp1} is the composition  of 
$ (L,B) \colon X \rightarrow Y\times Z $
with
$(L_\infty,B_\infty)\times Id  \colon  Y\times Z\rightarrow V\times W \times Z.$

Hence under the  assumptions of Theorem \ref{th:40}, $(L^1,B^1)$ is Fredholm, and the same holds   for $(L^2,B^2)$ which is  a
composition of  $( L_\infty , B)$ with $( L, B_\infty )\times \Id.$   In particular,  \eqref{comp1} and  \eqref{comp2}  are elliptic boundary value problems, being ellipticity  equivalent to the Fredholm property of the induced operator. 

The above two decompositions give: 
\be \label{indice1}\ba{c}  \Ind(L^1,B^1)=  \Ind(L_\infty, B_\infty)+ \Ind(L,B)\\ 
\Ind(L^2,B^2)=  \Ind(L, B_\infty)+ \Ind(L_\infty,B).\ea\end{equation} 

Since  $\cL_\la(x,D)=\cL_\infty(x,D)$ for $x $ close to $\Gamma,$  we have $\cB^1=\cB^2.$ Moreover, the principal symbols of $\cL_\infty$ and  $\cL_\la$ commute   and hence the principal symbols of $\cL_\infty \cL$ and $\cL \cL_\infty$ coincide.  It follows  that  $\cL_\infty \cL- \cL \cL_\infty$  is of order $-1.$ Thus  the  families $(L^1,B^1)$ and 
$( L^2,B^2)$ differ by a family of compact operators  and therefore 

\be \label{indice2} \Ind(L^1,B^1)=\Ind(L^2,B^2).\end{equation}

Since $(L_\infty, B_\infty),$ is a constant family of operators of index $0,$  from \eqref{indice1}, \eqref{indice2}  we obtain: 
 
\be \label{indice} d\left( \Ind(L,B)\right)= d\left( \Ind(L^1,B^1)\right)=d\left(\Ind(L, B_\infty)\right)+ d\left(\Ind(L_\infty,B)\right).\end{equation}
The  degrees on the right hand side have been computed  in \cite {[Pe]}  and  in Step $2.$  Indeed, $d\left(\Ind(L_\infty,B)\right) $ is given by  \eqref{degre}  and $d\left(\Ind(L, B_\infty)\right)$ is given  by the same formula  
 involving the interior multiplicity $\mu_i(\cL,\cB),$  \cite[(4.25)]{[Pe]}). 
 
In conclusion: 
\begin{equation}\label{degol} d\left(\Ind (L,B)\right) = \begin{cases} \mu(L,B)  &\hif q\equiv 0 \mod 8 \\ \frac{1}{2}\mu(L,B)  &\hif  q\equiv 4 \mod 8,\end{cases}\end{equation}
  and the last argument in Step 2 completes the proof of the theorem. \qed

\section{An example} 
In this section we will show how  to construct  families of elliptic differential boundary value problems verifying the hypotheses of Theorem \ref{th:40}. 
   Examples  with  pseudo-differential boundary conditions are  easy to find.   However, exhibiting  concrete examples with  differential boundary conditions is far from being simple.  Indeed, very little is known about the set of elliptic systems of differential operators of a given order, and even less about the structure of elliptic boundary values problems. 
   
    Following a suggestion of Atiyah in \cite{[At-9]}, we will take an indirect approach by  approximating the principal symbol   of a family of elliptic  pseudo-differential  operators with symbols of  families of elliptic differential operators of sufficiently high order.   Atiyah's  idea is to consider  the set of  elliptic symbols  $A(n,r,2k),$  whose elements  are $r\times r$ matrices with homogeneous polynomial entries of order $2k$ in variables $(\xi_1\dots,\xi_n).$ Then the approximations are constructed  using the fact  that the set of restrictions to the unit sphere of elements of $A(n,r) =  \cup_{k\geq 0} A(n,r,2k)$  is dense in the set of even continuos functions from the sphere into $GL(r; \C)$ (see also \cite{[Va]} for  a related result).
   
 \sk
 Before going to this point we need some preliminaries. 
 First, let us observe that  the reality condition  on the principal symbols is irrelevant  to the validity of  \eqref{cher2}. Hence  \eqref{cher2}  holds true not only for the map $\tau$ defined in \eqref{parfed} but  in general. Thus, for  any smooth  map $ \phi \colon S^q\times S(\gG) \ra GL(r,\C)$ with $q =2t$ even, we have: 
   \begin{equation}\label{cher4} 
  <p_*\ch([\phi]);[S^q]>= \deg{\phi}. 
  \end{equation} 
  
  Secondly, since we are dealing with homotopy invariants  of maps  with values  $GL(r;\C)$  and since the unitary group $U(r)$ is a deformation  retract of $GL(r,\C),$ in our  discussion we can safely assume that $\phi$ takes values in the unitary group $U(r).$

 The traces $\theta_i= tr(u^{-1}du)^{2i-1}$ of the  odd powers  of  the Maurer-Cartan  matrix-differential form   $u^{-1}du$  of $U(r)$ are bi-invariant and hence harmonic differential forms. The forms  $\theta_i$  define  cohomology classes  $[\theta_i]\in  H^{2i-1}(U(r);\C)$ in the de  Rham cohomology with coefficients in $\C$ which are  known to be  generators of  the exterior algebra $H^{odd}(U(r);\C).$  The pullback of $ \theta_{q+2n-3}$ by the map $\phi$ is  the form $\phi^{-1}d\phi,$ and hence   we can write 
  \begin{equation}\label{degint}   \deg(\phi)=N \displaystyle\int_{S^q\times S(\gG)} \phi^* \theta_{q+2n-3}. \end{equation}
   Here $\phi^*$ denotes  the pullback of $\phi$ and $N$  is as in \eqref{bottfed}.

Finally,  let us recall  that if $ \psi  \colon S^{2v-1} \rightarrow Gl(r;\C),$ is any continuous map, the clutching construction associates to $\psi$  a vector bundle $[\psi]$  over $S^{2v}$  obtained by gluing  via the map $\psi$   two trivial complex bundles of rank $r$  over the upper and lower  hemispheres $D_\pm$ of $S^{2v}.$   For  $r\geq v,$ the above  construction induces an isomorphism of  $\pi_{2v-1}\big(U(r)\big) $ with $\tilde K(S^{2v}).$

    In  3.1 of \cite{[Fe]} Fedosov showed  that for smooth  $\psi $ and an appropriate choice of orientation of $S^{2v-1}$  

 \begin{equation}\label{degeta} <\ch ([\psi]);[S^{2v}]>=\deg(\psi)=  \displaystyle N \int_{ S^{2v-1}} tr (\psi^{-1}d\psi)^{2v-1}.\end{equation} 
  
 \sk

With this said,  let us go to the example.
 
Let $n\geq 3 $ be odd. For simplicity,  choose $\Omega$ such  that the cotangent bundle of $\gG$ is trivial, e.g.,  take as $\Omega$ the region bounded by  an  $(n-1)$-torus  $\gG= (S^1)^{n-1}.$   Then $S(\gG)\cong \gG \times S^{n-2}.$ 
  Consider the map  $f$ defined as the composition
\begin{equation}\label{fdef}\displaystyle \left.\begin{array}{ccccc} & \Id\times \pi &  & g &  \\S^q\times \gG \times S^{n-2} & \to & S^q\times \gG \times \R P^{n-2} & \to & S^{q+2n-3}\end{array}\right.,\end{equation}
where $\R P^{n-2}$ is the real projective space, $\pi \colon S^{n-2} \ra \R P^{n-2}$ is  the canonical projection and  $g\colon S^q\times T^{n-1} \times \R P^{n-2} \ra S^{q+2n-3}$ is a smooth map having  Brouwer  degree one and  sending  $\{\infty\} \times \gG \times \R P^{n-2}$ into a point.  

Notice that  Brouwer's degree  $\deg_B\,g$ is  defined  because  $\R P^{n}$ is orientable for odd $n.$  Moreover,  $\deg_B\,f=2$  since $\deg_B\, \pi =2,$  in this case.
  
Choose an $r\geq q+2n-3$  and define   $\phi \colon S^q\times \gG \times S^{n-2} \ra GL(r;\C) $ to be  the composition of  the map $f$ with a  map  $\psi \colon S^{q+2n-3}\ra GL(r;\C) $ representing a generator of   $\pi_{q+2n-3}(U(r)) \simeq \tilde K(S^{q+2n-2}) \simeq \Z .$   By construction the map $\phi$ is even, in the  variable $\xi,$  i.e., $\phi(\la,x', -\xi')= \phi(\la,x', \xi').$  

By the Change of Variables Theorem,  

\begin{equation} \deg(\phi)=N \displaystyle\int_{S^q\times S(\gG)} \phi^* \theta_{q+2n-3} = 
  N\, \deg_B\, f\, \displaystyle\int_{S^{q+2n-3} } \psi^* \theta_{q+2n-3}= 2\deg(\psi)=2, \end{equation} 
being $\deg(\psi)=1$ by \eqref{degeta}.

\begin{lemma}\label{diffapprox} The map  $\phi\colon S^q\times \gG \times S^{n-2} \ra U(r)$
constructed  above,  considered as a map with values in $Gl(r;\C),$ can be uniformly approximated  by  the restriction to  $S^q\times \gG \times S^{n-2}$  of the symbol of a family of homogeneous elliptic differential operators.\end{lemma}

\proof  Let  $\cC= C(S^{n-2}; \C^{r\times r}),$ endowed with the sup norm, and let  $ \cA \subset \cC $ be the set of all  restrictions to $S^{n-2}$ of polynomial maps from $\R^{n-1}$ to $\C^{r\times r}.$  By  the Stone-Weierstrass theorem $\cA$   is dense in $\cC.$  Using this and  smooth  partitions of unity  on $S^q\times \gG,$  for every $\epsilon>0$ we can find a smooth  map 
$$\rho \colon S^q\times \gG \times \R^{n-1} \ra \C^{r\times r} $$ such that  \begin{equation} \label{approxrho} \| \rho(\la,x',\xi')- \phi(\la,x', \xi')\|_\infty <\epsilon \,\hforall\, (\la,x', \xi')\in S^q\times \gG \times S^{n-2} \end{equation}
and such that  $\rho(\la,x', \xi') = \sum_{|\alpha|\leq t} a_\alpha(\la,x') (\xi')^\alpha.$

 For simplicity let us assume that $t=2s$ is even. We rewrite the last expression in the form  $\rho(\la,x', \xi') =\sum_{i=0}^{2s}  h_i(\la,x', \xi'),$  where $h_i(\la,x', \xi') = \sum_{|\alpha|= i} a_\alpha(\la,x') (\xi')^\alpha.$ Thus the maps $h_i$ are homogeneous polynomials in $\xi$  of degree $i.$ 
 
 Since $\phi(\la,x',\xi') = \frac{1}{2}\,[ \phi(\la,x', \xi') +\phi(\la,x',- \xi') ],$  the restriction to  $S^q\times\gG \times S^{n-2}$ of the   even part     $\rho_{ev} (\la,x', \xi')= \sum_{i=0}^ {s}  h_{2i}(\la,x', \xi')$ of $\rho$ also verifies \eqref{approxrho}. Now  we can approximate   $\phi$  by the restriction of a map $h$  which is a homogeneous polynomial in $\xi$.   In fact,  at points with $|\xi|=1 $  the values of the even homogeneous polynomial map $ h(\la,x', \xi')=  \sum_{i=0}^{s}   |\xi |^ {2i} h_{2i}(\la,x', \xi'),$  coincide with those of  $\rho_{ev}.$ Thus  $h$ is the  symbol of  a family of homogeneous differential operators which becomes elliptic after choosing an $\epsilon$ small enough.   \qed 

\begin{remark} {\rm  In the above lemma  the triviality of the  cotangent bundle is inessential.  Indeed,  the  proof shows that any parametrized family of even   maps  from the  cotangent sphere bundle to $GL(r;\C)$  can be uniformly approximated by a family of principal symbols of elliptic differential operators.} \end{remark} 

  Choose $m,l$   such that $ml=r$  and consider the boundary value problem $(\cL^0,\cB^0),$ where  $\cL^0=( \Delta^{l}+ \mu) \Id_m,$ acting on $\C^m$ valued functions  and $\cB^0=(\ga_0,...,\ga_{l-1})$ is  the  Dirichlet boundary condition of order $l-1.$
  
   Since  $\Delta^{l}$   is a  strongly elliptic operator, taking $\mu $ big enough, we can assume that  the operator $(L^0,B^0)$  induced  by $(\cL^0,\cB^0)$ on Hardy-Sobolev spaces is an isomorphism.  Let us denote with $\cH=\{\cH_\la : \la \in S^q\}
  $  the family of  homogeneous elliptic differential operators on $\gG$ associated to the symbol $h$ constructed in  Lemma \ref{diffapprox}. 
  
   Since the index of an operator depends only on the homotopy class of the restriction of its principal symbol to the unit sphere bundle $S(\gG),$ and since the restriction of $h_\infty$  can be taken arbitrarily close to the  constant symbol $\phi_\infty$ we have that $\ind H_\infty=0.$   By eventually taking a lower order perturbation, we can also assume that the operator $H_\infty$ induced by $\cH_\infty$ is an isomorphism. 
  
 Let us consider now  the family  $(\cL,\cB)$ with  $\cL= \cL^0$ constant and $\cB=\cH\circ \cB^0.$  The Fredholm property of the induced operator  is equivalent to the ellipticity of the boundary value problem. Hence writing  $(L,B)$   in the form  \[(L,B) =( \Id \times H)\circ(L^0,B^0)\]   we see that  the family $(\cL, \cB) $ is a family of elliptic differential boundary value problems with complex (matrix) coefficients.   Moreover we have that $(L_ \infty, B_\infty)$ is an isomorphism.   As a family of complex differential operators of index $0,$ taking $k=2l$ the induced family   $(L,B)\colon S^q \times \hk \ra \hs\times\hbb$  has an index bundle $\Ind (L,B) \in \tilde K(S^q).$
 
 We have  
 \be \label{indice5} \Ind(L,B)=  \Ind(\Id\times H)+ \Ind(L^0,B^0)=\Ind H,\end{equation} being $(L^0,B^0)$ constant. 
 
 Since  the restriction of $h$ to $S^q\times S(\gG)$ homotopic to $\phi,$  from \eqref{indice5} we obtain 
\be \label{indice6} <\ch(\Ind(L,B));S^q>=<\ch(\Ind H);S^q>= <p_*\ch(\phi);S^q>= \deg(\phi) = 2.  \end{equation} 

 Now, we identify $\C^m$ and $\C^r$ with $\R^{2m}$ and $\R^{2r}$ respectively and consider $(\cL,\cB)$ as a family of real differential operators. In order to avoid confusions we will denote this family  with $(\cL^r,\cB^r)$. The ellipticity of this family is a consequence of the ellipticity of the corresponding complex family.  By construction the operator    $(\cL^r_\infty,\cB^r_\infty )$ verifies $H_2,$  and being $\cL^r$ constant,  the restriction of this family to parameters belonging to  $\R^q$ verifies  the assumptions $H_1-H_3.$
 
 Our aim is to apply Corollary \ref{lowerorder}   to the (real) nonlinear perturbations of the restricted family.  For this we have  to evaluate the multiplicity $\mu_b(\cL^r,\cB^r).$ 
 
  Much as in the proof of Theorem \ref{th:41}  we compute it from  the degree of the index bundle  of the complexification $ (L^{rc},B^{rc})$  of  $(L^r,B^r).$ 
  
Denoting  with $c\colon \tilde{KO} (-) \ra \tilde K(-)$  and $r \colon  \tilde K(-) \ra \tilde{KO} (-) $ the complexification and the realization homomorphism respectively we have:

\be \label{indice7} \mu_b(\cL^r,\cB^r) = <\ch(\Ind(L^{rc},B^{rc}));S^q>=  <\ch\big(c\circ r\, \Ind(L,B)\big);S^q>.\end{equation}  

The right hand side of \eqref{indice7} can be easily related to the left hand side of \eqref{indice6}. 
Indeed, for $q=4s,$ by  \cite[Theorem 30]{[Mi]}, the Chern classes of $c\circ r(\eta)  \in \tilde K (S^q)$  verify  $c_{2s}(c\circ r\,(\eta)) = \pm 2c_{2s}(\eta).$ 
Using  $c_{2s}(\eta)= \pm (2s-1)! \ch(\eta)$ we conclude from \eqref{indice6}  that 
 
\be \label{indice8}  \mu_b(\cL^r,\cB^r)= \pm 2<\ch\big(\Ind(L,B)\big);S^q>= \pm 2 \deg(\phi) =\pm 4.\ee

The function $n(q)$ defined in \eqref{muca}  always assumes values greater or equal than $24.$ Hence, by Corollary \ref{lowerorder}, 
any family   \[(\cf,\cg) \colon \R^q \times C^\infty( \bar\Omega;\R^{2m}) \rightarrow C^\infty( \bar\Omega;\R^{2m}) \times\ C^\infty(\partial\Omega;\R^{2r})\]
of the form 
 \[ (\cL^r_\la(x,D) u + \cf'(\la, x,u,\ldots), \cB^r_\la(x,D)u + {\cg'}(\la, x,u,\ldots)), \]  where $(\cf',\cg')$ is any lower order perturbation with  $\cf'(\la, x,0)=0,\cg'(\la, x,0)=0 $ and  such that the coefficients of the linearization of $(\cf',\cg')$ converge  uniformly  to $0$ as $\la\ra \infty,$ will have at least one  bifurcation point. 
 \sk
In conclusion, the above construction of linear elliptic  boundary value problems with  multiplicity smaller in absolute value than any $n(q)$  produces  examples to which our bifurcation criteria can be applied for all data  $(q,n,m,l)$  with $n$ odd, $q$ divisible by $4$ and   $ml=r \geq q+2n -3.$ The same method can be used in order to construct families  of linear differential operators with constant boundary conditions  and $\mu_i(\cL,\cB)$ small but nonzero.    

\section{Appendix}

   We are going to discuss  an alternative description of $\deg(\phi)$   taken from  Atiyah and Singer \cite[Section 9]{[At-Si-III]} which works also for maps $\phi$  that are only continuous and compare  it  with Fedosov's approach used in this paper.   This  sheds some light  on the  construction of families of invertible matrices with a  given degree.
       
  The left  hand side  of  \eqref{cher4} can be computed in a different way.  This is done  in \cite[section 9]{[At-Si-III]}  using  generators of $H^*(U(r);\Z)$  which transgress to  the universal Chern classes in $H^*(BU(r)).$  Without going into details, which, by the way, are similar to the proof of Theorem  \ref{th:41}, the  conclusions  are as follows: 
 
As before, let  $V= S^q\times S(\gG)$ and $v= q/2 +n -1.$  Using $U(r)/U(r-1)=S^{2r-1}$ it is easy to see  that, if $r\geq v,$  every continuous map $\phi \colon V \rightarrow U(r)$ is homotopic to one of the form $\diag(  \phi', \Id_{r-2v+1}),$ where $\phi'$  takes its  values in $U(v).$  Taking the first column   $\phi'_1$ of $\phi'$ we obtain a map $ \phi'_1  \colon V \rightarrow S^{2v-1}$ which,  being a map of oriented manifolds of the same dimension, has a well defined Brouwer's degree $ \deg_B\, \phi'_1\in \Z.$ It turns out that   $\deg_B\, \phi'_1$  is  divisible by $(v-1)!.$

Define 
 \begin{equation} \label{bottdeg'} \deg'(\phi)=  \displaystyle \begin{cases}  \frac{1}{(v-1)!}\, \deg_B\, \phi'_1 & \, \hif  \, r\geq v\\ 0 &\,\hif \,r< v. \end{cases}\end{equation} 
 
%\begin{remark} {\rm In the case  $V=S^{2v-1}$   \eqref{bottdeg'}  is called  Bott's degree  in \cite {[At-8]} because  the assignment   $\phi \ra \deg'(\phi)$  induces  the Bott isomorphism  $ \pi_{2v-1}\big(U(r)\big)\simeq \Z.$}\end{remark} 

It  is shown in \cite[corollary(9.5) ]{[At-Si-III]} that   
\begin{equation} \label{indias}  \deg'(\phi)=<\ch ( [\phi]);[S^q\times S(\gG)]> \end{equation}  (the statement  in Corollary $ (9.5)$ is in the case without  parameters  but the proof with parameters is the same). 

Since  the direct image homomorphism $f_*$ commutes with composition of maps  and since $cte_* $ coincides with the evaluation on the fundamental class  one easily verifies that 
\begin{equation} \label{?} <\ch( [\phi]);[S^q\times S(\gG)]>  =<p_* \ch [\phi];[S^q]>. \end{equation} 

The above relation, together with \eqref{cher4} and \eqref{indias},  allows us  to conclude that  if $\phi$ is smooth,   then $ \deg'(\phi)= \deg(\phi).$  In particular  $\deg'(\phi)$   extends   $\deg(\phi)$ to all continuous maps from $S^q\times S(\gG)$ to $GL(r;\C).$ 
\sk 
 As a side remark let us observe that the use of to above extension of the degree applied to the map $\sigma$  allows us relax assumption $H_3$ in Theorem \ref{th:40} to\,: 
\sk
$H'_3)$ - The restrictions to $ \partial \Omega$    of  the coefficients   of the  leading terms  of  ${\mathcal L}_\la (x,D)$
 are independent of $\la.$
\sk

Indeed, we  have only  to improve  Step $1$ of the proof of  Theorem \ref{th:40}.  Under  assumption  $H'_3,$ the construction  of  $\sigma$ in \eqref{symbo} gives only a continuous map.  However,  using  $$
<\ch\Ind(L^c,B^c) ;[S^q]> = <p_* \ch [\sigma];[S^q]> = \deg'( \sigma) $$ and arguing as in  the  proof of Theorem $1.4.1$ in  \cite {[Pe]}  we  can show that $\J(\Ind L)\neq 0$ whenever $\deg'(\sigma)$ is not divisible by $n(q).$

  On the negative side,  let us observe that, while  $\deg$ is given explicitly by an integral of a differential form, the definition of $\deg'$  in \eqref{bottdeg'} is far less explicit. Notice also that we could define $\deg'$ using approximations of continuous maps by the smooth ones,  as is often  done for the Brouwer degree.

\end{document}